\documentclass[a4paper,12pt]{amsart}

\usepackage[T1]{fontenc}

\usepackage{amssymb}
\usepackage{amsmath}

\usepackage{enumitem}
\newtheorem{theorem}{Theorem}
\newtheorem{corollary}[theorem]{Corollary}
\newtheorem{lemma}[theorem]{Lemma}
\newtheorem{proposition}[theorem]{Proposition}

\theoremstyle{remark}
\newtheorem{remark}{Remark}

\newcommand{\Z}{{\mathbb{Z}}}

\begin{document}
%\linenumbers
\title[On the arithmetic of polynomials]{On the arithmetic of polynomials with coefficients in  Mordell-Weil type groups.}

\author{Stefan Bara\'{n}czuk} 
\address{ Faculty of Mathematics and Computer Science, Adam Mickiewicz University, ul. Umultowska 87, Pozna\'{n}, Poland}
\email{stefbar@amu.edu.pl}
\begin{abstract}In this paper we prove  Hasse local-global principle for polynomials with coefficients in Mordell-Weil type groups over number fields like $S$-units, abelian varieties with trivial ring of endomorphisms and odd algebraic $K$-theory groups.
\end{abstract}

\maketitle

\section{Introduction. Mordell-Weil type groups.}

We say that Hasse principle holds for $f \in \Z [x]$ when the following statement is true:
\begin{center}
	\textit{$f$ has an integer root if and only if it has a root over $\Z/p$ for almost all prime numbers $p$. }
\end{center} 
Hasse (\cite{H2}) has proved, correcting a statement by Wegner, that if $f \in \Z [x]$ is irreducible and has a root over $\Z/p$ for almost all prime numbers $p$ then it has a rational root (later Fujiwara has obtained this result for some families of reducible polynomials, see \S 1 of \cite{F}) thus the principle holds for irreducible monic polynomials with integer coefficients.\\

%In particular it means that a polynomial of degree $2$ or $3$ with integer coefficients  has an integer root if and only if it has a root modulo all but finitely many prime numbers. Later Fujiwara (\cite{F}, Proposition 2) extended the result to degree $4$ polynomials (see also \cite{BB}, Remark 2). 

The main aim of this paper is to extend Hasse principle to polynomials with coefficients in  Mordell-Weil type groups (Theorem \ref{main}). As a by-product of the proof, we obtain a dynamical  principle, answering a question by Silverman (Corollary \ref{Silv}). \\ 
 
For the purpose of this paper we define Mordell-Weil type groups via the following abstract nonsense axiomatic setup:

Let $B$ be an abelian group and $r_{v} \colon B \to B_{v}$ be family of groups homomorphisms indexed by the set of primes $v$ in a number field $\mathbb{K}$ whose targets $B_{v}$ are finite abelian groups. We will use the following notation:\\
\begin{tabular}{ll}
$P \mod v$	&	denotes $r_{v}(P)$ for $P \in B$\\
$P=Q \mod v$ & means  $r_{v}(P)=r_{v}(Q)$ for $P,Q \in B$\\
$\Lambda_{\mathrm{tors}}$ &   the torsion part of a subgroup $\Lambda < B$ \\
$\mathrm{ord} \ T$ & the order of a torsion point $T \in B$ \\
$\mathrm{ord}_{v} P$ &  the order of a point $P \mod v$ \\
$l^{k} \parallel n$	& means that $l^{k}$ exactly divides $n$, i.e. $l^{k} \mid  n$		and	$l^{k+1} \nmid n$ \\ & where  $l$ is a prime number, $k$ a positive integer \\ & and $n$ a natural number.						
\end{tabular}\\
We impose the following two assumptions on the family $r_{v} \colon B \to B_{v}$:

\begin{enumerate}[label=A{\arabic*}]
	\item \label{o posylaniu} Let $l$ be a prime number and $(k_{1}, \ldots, k_{m})$ a sequence of nonnegative integers. If $P_{1}, \ldots, P_{m} \in B$ are  points linearly independent over $\Z$ then there is a positive density set of primes $v$ in $\mathbb{K}$ such that $l^{k_{i}} \parallel \mathrm{ord}_{v} P_{i}$ if $k_{i}>0$ and  $l \nmid \mathrm{ord}_{v} P_{i}$ if $k_{i}=0$.
	\item \label{o torsyjnych} For almost all $v$ the map $B_{\mathrm{tors}} \to B_{v}$ is injective.
\end{enumerate}

This axiomatic setup holds for the following families of groups (see e.g. \cite{Bar2}):
\begin{itemize}
	\item ${R}_{\mathbb{K},S}^{\times}$, $S$-units groups, where $\mathbb{K}$ is a number field and $S$ is a finite set of ideals in the ring of integers ${R}_{\mathbb{K}}$,
	\item $A(\mathbb{K})$, Mordell-Weil groups of abelian varieties over number fields $\mathbb{K}$ with $\mathrm{End}_{\bar{\mathbb{K}}} (A) = \mathbb{Z}$,
	\item $K_{2n+1}(\mathbb{K})$, $n>0$, odd algebraic $K$-theory groups.
\end{itemize}	
\section{Main results.}

Before formulating our main theorem we define the greatest common divisor of a polynomial with coefficients in  Mordell-Weil type groups. 
Let 
\[F(n)=P_{0} + n P_{1}+ \ldots +n^{d} P_{d}\]
be polynomial with coefficients $P_{0}, P_{1},\ldots, P_{d}\in B$. Denote the order of the group $B_{\mathrm{tors}}$ by $t$ and write 
\[ t F(n)=f_{1}(n)G_{1}+\ldots +f_{i}(n)G_{i}\]
for linearly independent points $G_{1}, \ldots, G_{i} \in B$ and nonzero polynomials $f_{1},\ldots, f_{i}\in \Z [n]$. We define 
\textit{the greatest common divisor of} $F$ to be  $\gcd (f_{1},\ldots, f_{i})$ when $t F \ne 0$ and $0$ when $t F = 0$.

We say that Hasse principle holds for $F$ when the following conditions are equivalent:
\textit{\begin{enumerate}[label=C{\arabic*}]
		\item \label{cloc} For almost every $v$ there exist a rational integer $n$  and  $T \in B_{\mathrm{tors}}$ such that
		\begin{equation}\label{zal}
		F(n) = T \mod v.
		\end{equation}
		\item \label{cglob}
		There exist a rational integer $n$ and  $T \in B_{\mathrm{tors}}$ such that 
		\[	F(n) = T.\]
	\end{enumerate}}
	
\begin{theorem}\label{main}
	 Hasse principle (as presented above) holds for  $F$ if and only if Hasse principle (as presented in the Introduction) holds for the greatest common divisor of $F$.
\end{theorem}

\proof
 Every $r_{v}$ is homomorphism so condition \ref{cglob} implies \ref{cloc} thus it remains to investigate the converse implication.\\
 
If all $P_{0}, P_{1},\ldots, P_{d}$ are torsion then the result follows immediately from Assumption \ref{o torsyjnych}. So suppose one of them is nontorsion.
 
Following the notation in the above definition of the greatest common divisor we multiply \eqref{zal} by $t$ and get that for almost every $v$ there exists a rational integer $n$  such that
		\begin{equation}\label{zalnontor}
		f_{1}(n)G_{1}+\ldots +f_{i}(n)G_{i}=0 \mod v.
		\end{equation}
We have to analyse when there exists a rational integer $n$  such that
\[f_{1}(n)G_{1}+\ldots +f_{i}(n)G_{i}=0.\]	

Let us first consider the case $i=1$. Rewrite \eqref{zalnontor} as		
\begin{equation}\label{zalnontor1}
		f(n)G=0 \mod v.
		\end{equation}		
By Assumption \ref{o posylaniu} for every prime number $l$ there are infinitely many $v$'s such that $l \mid \mathrm{ord}_{v} G$ so by \eqref{zalnontor1} there exists a rational integer $n$ such that $l \mid f(n)$ thus $f$ has an integer root if and only if  Hasse principle holds for $f$.\\

For the proof in the case $i>1$	rewrite \eqref{zalnontor} as
	\begin{equation}\label{zalnontormany}
		\gcd(f_{1},\ldots, f_{i})(n)\big( g_{1}(n)G_{1}+\ldots +g_{i}(n)G_{i}\big) =0 \mod v.
		\end{equation}
Applying Lemma \ref{lemat} to polynomials $g_{1}, \ldots, g_{i}$ we can assume without the loss of generality that for almost every prime number $l$ there exist  rational integers $m_{2}, \ldots, m_{i}$ such that the polynomial 	
\[h=g_{1}+m_{2}g_{2}+\ldots+m_{i}g_{i}\]
has no root over $\Z/l$. Write $ g_{1}(n)G_{1}+\ldots +g_{i}(n)G_{i}$ in \eqref{zalnontormany} as
\[h(n)G_{1}+g_{2}(n)(G_{2}-m_{2}G_{1})+\ldots+g_{i}(n)(G_{i}-m_{i}G_{1}).\]
The points $G_{1},\ldots, G_{i}$ are linearly independent so are the points \[G_{1}, G_{2}-m_{2}G_{1}, \ldots, G_{i}-m_{i}G_{1}\] and by Assumption \ref{o posylaniu} there are infinitely many $v$'s such that 
\[\begin{array}{ccccc}
l \mid \mathrm{ord}_{v}G_{1} &  \mathrm{and} & l \nmid \mathrm{ord}_{v}(G_{j}-m_{j}G_{1}) &  \mathrm{for} & j \in \left\lbrace 2, \ldots , i \right\rbrace 
\end{array}\]
so \[l \mid \mathrm{ord}_{v}\big(g_{1}(n)G_{1}+\ldots +g_{i}(n)G_{i}\big)\]
thus by \eqref{zalnontormany} there exists a rational integer $n$ such that $l \mid \gcd(f_{1},\ldots, f_{i})(n)$ hence  $\gcd(f_{1},\ldots, f_{i})$ has an integer root
if and only if  Hasse principle holds for $\gcd(f_{1},\ldots, f_{i})$.	
\qed

\begin{remark}\label{deg gt 4}
By the result of Hasse and  Proposition 2 in \cite{F} Hasse principle holds for all polynomials with coefficients in Mordell-Weil type groups whose greatest common divisor is a monic polynomial of degree $\le 4$. This is no longer true if we increase the degree or drop the assumption of monicity.
	
	The polynomials 
	\[f ( n )=( n^3 - 19)( n^2 + n +1)
	\]
	and 
	\[f ( n )=( n^2- 13)( n^2 - 17)( n^2 - 221)
	\]
	(\cite{BB}, Examples 1 and 2) have no integer roots but have roots modulo every integer and multiplying them by an appropriate power of $n^2+1$ we extend these examples to polynomials of arbitrary degree $>4$. 
		
		For every $d\ge 2$ the polynomial
		\[f(n)=(3n-2)(2n-3)^{d-1}\]
		has no integer root but has a root modulo every integer.

\end{remark}
\begin{remark}\label{T eq 0}
In view of the original Hasse principle one may ask if it is possible to fix a torsion point in Theorem \ref{main} when the group $B_{\mathrm{tors}}$ is nontrivial, i.e., whether we have equivalence of the following statements:
\begin{itemize}
		\item \textit{For almost every $v$ there exists a rational integer $n$   such that} 
		\[
		F(n) = 0 \mod v.
		\]
		\item \textit{There exists a rational integer $n$ 
		 such that} 
		\[
			F(n) = 0.
		\]
	\end{itemize}
Analysing the proof of Theorem \ref{main} we see that the answer is positive for polynomials of the form
\[F(n)=f_{1}(n)G_{1}+\ldots +f_{i}(n)G_{i}\]
where $G_{1}, \ldots, G_{i}$ are linearly independent points such that Hasse principle holds for the greatest common divisor of $F$. 

However, the answer is negative in general. In fact for every group $B$ with a nontorsion point $P$ and nontrivial torsion subgroup we can construct a counterexample of arbitrary degree $d \ge 2$. Indeed, let $P\in B$ be nontrosion and $T\in B$ be torsion with $\mathrm{ord} \ T >1$.  Consider the polynomial 
\[F(n)=n^{d-1}\left( \left(  \left(\mathrm{ord} \ T + 1 \right) n -1 \right) P   +T\right).  \]  
For every rational integer $n$ we have $F(n) \ne T$  but for every $v$ there exists a rational integer $n$ such that 
\[
		F(n) = T \mod v.
		\]
Indeed, for every $v$ fix natural numbers $k,l$ such that $\gcd (l, \mathrm{ord} \ T)=1$ and  $\mathrm{ord}_{v} P$ divides $l (\mathrm{ord} \ T)^{k} $. Our $n$ is a solution of the following system of congruences 
\[\left\{ \begin{array}{l}
\left(\mathrm{ord} \ T + 1 \right) n -1 \equiv 0 \mod (\mathrm{ord} \ T)^{k} \\
n \equiv 0 \mod l .
\end{array} \right.
\]
Observe that $F(0)=0$ so this example does not collide with Theorem \ref{main}.		\\

As for the case of degree $1$ polynomials, note that the condition \[P+nQ=0\] equals \[P \in \left\langle Q\right\rangle\]
where $\left\langle Q\right\rangle $ denotes the subgroup of $B$ generated by $Q$. We treat this case in the next Remark.
\end{remark}
\begin{remark}
Theorem \ref{main} might be viewed as a variation of principles known as \textit{detecting linear dependence} addressed recently in numerous papers; we do not discuss them here and refer to \cite{Bar2} instead. Let us shortly say that they are local-global principles for the statement \[P\in \Lambda \]where $P\in B$ and $\Lambda$ is a subgroup of $B$. In all  results of this kind  $P$ was assumed to be nontorsion; however in regard to the counterexample constructed in the previous remark we should note that this assumption is not necessary as we demonstrate in the  Proposition below. 
\end{remark}
\begin{proposition}
Let $P$ be a point in $B$ and $\Lambda$ a subgroup of $B$. Then $P \in \Lambda$ if and only if $P \in \Lambda$ modulo almost all $v$. 
\end{proposition}
\proof
	The statement is known for $P$ nontorsion (see the Theorem in \cite{Bar2}) thus it remains to prove it for $P$ torsion.\\
	($\Rightarrow$) Every $r_{v}$ is homomorphism so if $P \in \Lambda$ then $P \in \Lambda$ modulo $v$.\\ 
	 ($\Leftarrow$) Let $k$ be the least natural number such that $k P \in \Lambda_{\mathrm{tors}}$. If $k=1$ we are done. So suppose $k \ge 2$. For every $T \in \Lambda_{\mathrm{tors}}$ we have $k \mid \mathrm{ord}(T-P)$. Fix a prime number $l$ dividing $k$. Let $G_{1},\ldots,G_{h}\in \Lambda$ be linearly independent points that together with points in $\Lambda_{\mathrm{tors}}$ generate $\Lambda$. By Assumption \ref{o posylaniu} there is a positive density set of primes $v$ such that $l$ does not divide $\mathrm{ord}_{v}G_{i}$ for $i=1,\ldots,h$ hence by Assumption \ref{o torsyjnych} it follows that $P \notin \Lambda$ modulo $v$ for almost all such $v$'s. 
\qed

\section{A dynamical result.}

In regard to   Proposition 4 of \cite{Bar3}, Joseph Silverman asked about local-global principle for sequences of the form
\[S_{n}=\varphi^{n}F(n)\]
where $\varphi$ is an endomorphism (i.e., multiplication by an integer in case of Mordell-Weil type groups we deal with in this paper) and $F$ a polynomial with coefficients in Mordell-Weil type groups. We answer to the question with the following Corollary.
\begin{corollary}\label{Silv}
The following are equivalent if and only if Hasse principle holds for the greatest common divisor of $F$:	
	\begin{enumerate}[label=C{\arabic*}]
		\item  \label{c1loc} For almost every $v$ there exist a natural number $n$  and  $T \in B_{\mathrm{tors}}$ such that
		\begin{equation*}
		\varphi^{n}F(n) = T \mod v.
		\end{equation*}
		\item \label{c1glob}
		There exist a rational integer $a$ and  $T \in B_{\mathrm{tors}}$ such that 
		\begin{equation*}
			F(a) = T.
		\end{equation*}
	\end{enumerate}	
\end{corollary} 

\proof Observe that condition \ref{c1glob} always implies condition \ref{c1loc}. Indeed, let $a$ be a rational integer such that $F(a) = T$ for some  $T \in B_{\mathrm{tors}}$. Denote the exponent of the group $B \mod v$ by $e_{v}$. Then we can take $n$ to be $a \mod e_{v}$. \\
So we are reduced to investigation of the implication \ref{c1loc} $\Rightarrow$ \ref{c1glob}. But this we get rewriting the proof of Theorem \ref{main} omitting $l$'s dividing $\varphi$. 
\qed
\begin{remark}
Remarks \ref{deg gt 4} and \ref{T eq 0} remain valid if we modify their statements as in Corollary \ref{Silv}. In particular, taking $\varphi$ such that $\varphi$ equals $1$  modulo the order of $B_{\mathrm{tors}}$ we get from Remark \ref{T eq 0} examples showing that the condition
\begin{center}
\textit{for almost every $v$ there is a natural number $n$ such that $S_{n}=0 \mod v$} 
\end{center}
does not have to imply that $S_{n}=0$ for some natural number $n$. 
\end{remark}

\section{Auxiliary results.}

\begin{proposition}\label{stw}
Let  $p$ be a  prime number and  $f,g\in \mathbb{F}_{p}  [x]$ polynomials having  no common root.  There exist $a,b\in \mathbb{F}_{p}$	such that the polynomial
\[h=af+bg\]
has no root.
\end{proposition}
\proof
If one of $f,g$ has no root we are done. So suppose that both $f,g$ have a root. Since by our assumption $f,g$ have no common root the set 
\[S=\left\lbrace x\in \mathbb{F}_{p} \, \colon f(x)\ne 0 \land g(x)\ne 0  \right\rbrace \]
has at most $p-2$ elements. If $S=\emptyset$ then put $h=f+g$.  If $S\ne\emptyset$ consider the set 
\[Q=\left\lbrace {f(x)}/{g(x)} \, \colon x \in S\right\rbrace. \]
Choose $n\in \mathbb{F}_{p} ^{\times} \setminus Q$ and put $h=f-ng$.
\qed

\begin{lemma}\label{lemat}
	For $f_{1}, \ldots, f_{k} \in \Z[x]$, $k>1$ let $\gcd (f_{1}, \ldots, f_{k})=1$. Then for almost every prime number $p$ there exist $c_{1}, \ldots, c_{k} \in \Z$ such that polynomial
	\[c_{1}f_{1}+\ldots +c_{k}f_{k}\]
	has no root over $\mathbb{F}_{p}$. 
\end{lemma}
\proof
	The proof is by induction on $k$.\\
	
Basis: If $f_{1},  f_{2} \in \Z[x]$ are coprime then  for almost every prime number $p$ they have no common root over $\mathbb{F}_{p}$ and we apply Proposition \ref{stw}.\\

Inductive step: Suppose that the Lemma is true for $k-1$. Denote $w=\gcd(f_{1}, \ldots, f_{k-1})$ and for every $i\in \left\lbrace 1, \ldots, k-1\right\rbrace $ write $f_{i}=w {g_{i}}$. Since $\gcd( g_{1}, \ldots,  g_{k-1})=1$ then by the assumption for almost every prime number $p$ there exist $a_{1}, \ldots, a_{k-1} \in \mathbb{F}_{p}$ such that the polynomial \[h=a_{1}g_{1}+ \ldots+ a_{k-1}g_{k-1}\] has no root in $\mathbb{F}_{p}$. Since $w$ and $f_{k}$ are coprime then  for almost every prime number $p$ they have no common root over $\mathbb{F}_{p}$ thus neither the polynomials $w h$ and $f_{k}$
 have it so by Proposition \ref{stw} there exist $s,t \in \mathbb{F}_{p}$ such that 
\[s w h + t f_{k} = s a_{1} f_{1}+\ldots+s a_{k-1} f_{k-1}+t f_{k}\]
has no root in $\mathbb{F}_{p}$.
\qed

\section*{Acknowledgements} The author is grateful to Dorota Blinkiewicz and Bart{\l}omiej Bzd\k{e}ga for stimulating discussions.

%-----------------BIBLIOGRAPHY----------

\bibliographystyle{plain}

\end{document}